\documentclass[12pt]{amsart}
\usepackage{graphicx} 
\usepackage{amsfonts, amsmath, amssymb, amsthm}
\usepackage{xcolor}
\usepackage{mathrsfs}
\usepackage{tcolorbox}
\usepackage{hyperref}
\usepackage{stix}
\usepackage{soul}
\usepackage[headheight=125pt,top=1.4in, bottom=1.4in, left=.8in, right=.8in]{geometry}
\title[Characterization of group Fourier transform]{Characterization of the group Fourier transform of Hilbert space Valued functions}

\author[Achar, Biswas, and Manna]{Satyajyoti Achar, Md Hasan Ali Biswas and Ramesh Manna}

\address{School of Mathematical Sciences, National Institute of Science Education and
Research, Bhubaneswar, An OCC of Homi Bhabha National Institute, Jatni 752050, India.}
\address{Homi Bhabha National Institute, Training School Complex, Anushakti Nagar, Mumbai 400094, India}
\email{satyajyoti.achar@niser.ac.in} 
\address{Department of Mathematics, Indian Institute of Technology
Gandhinagar, Palaj, Gandhinagar, Gujarat 382355, India}
\email{hasan.biswas@iitgn.ac.in}
\address{School of Mathematical Sciences, National Institute of Science Education and
Research, Bhubaneswar, An OCC of Homi Bhabha National Institute, Jatni 752050, India.}
\address{Homi Bhabha National Institute, Training School Complex, Anushakti Nagar, Mumbai 400094, India}
\email{rameshmanna@niser.ac.in}


\date{\today}

\begin{document}
\renewcommand{\Re}{\operatorname{Re}}
\renewcommand{\Im}{\operatorname{Im}}
\newcommand{\tr}{\operatorname{tr}}
\newcommand{\R}{\mathbb{R}}
\newcommand{\Hsn}{\mathbb{H}}
\newcommand{\Sl}{\mathcal{L}}
\newcommand{\A}{\mathcal{A}}
\newcommand{\F}{\mathscr{F}}
\newcommand{\N}{\mathbb{N}}
\newcommand{\Hi}{\mathcal{H}}
\newcommand*{\bigchi}{\mbox{\Large$\chi$}}
\newcommand{\He}{\mathbb{H}^n}
\newcommand{\blp}{\big (}
\newcommand{\brp}{\big )}
\newcommand{\Z}{\mathbb{Z}}
\newcommand{\Za}{\mathscr{Z}}
\newcommand{\T}{\mathbb{T}}
\newcommand{\J}{\mathbf{J}}
\newcommand{\C}{\mathbb{C}}
\newcommand{\diff}{\mathrm{d}}
\newcommand{\Df}{{\mathcal{F}_T}}
\newcommand{\D}{\mathcal{D}}
\newcommand{\B}{\mathcal{B}}
\newcommand{\M}{\mathcal{M}}
\newcommand{\G}{\mathcal{G}}
\newcommand{\nc}{\circledast_A}
\newcommand{\tc}{\otimes_\theta^\lambda}
\newcommand{\im}{\mathrm{Im}}
\newcommand{\re}{\mathrm{Re}}
\newcommand{\ch}{\mathcal{C_\theta}}
\newcommand{\nch}{\mathcal{C_{-\theta}}}
\newcommand{\f}{\mathcal{F}}
\def\Msp{{\boldsymbol M}}
\newtheorem{theorem}{Theorem}[section]
\newtheorem{lemma}[theorem]{Lemma}

\theoremstyle{definition}
\newtheorem{definition}[theorem]{Definition}
\newtheorem{example}[theorem]{Example}
\newtheorem{OP}[theorem]{Open problem}
\newtheorem{AP}[theorem]{Appendix}
\newtheorem{corollary}[theorem]{Corollary}
\newtheorem{xca}[theorem]{Exercise}
\newtheorem{propn}[theorem]{Proposition}
\theoremstyle{remark}
\newtheorem{remark}[theorem]{Remark}

\numberwithin{equation}{section}
\begin{abstract}
A remarkable theorem of S. Alesker, S. Artstein-Avidan, and V. Milman characterizes the Fourier transform on $\mathbb{R}^n$ as, in essence, the unique transformation on the space of tempered distributions that exchanges convolution with pointwise multiplication.
In this paper, we extend this perspective to the setting of the Heisenberg group $\He$, focusing on functions taking values in a separable Hilbert space $\Hi$. We establish a characterization of the group Fourier transform on $\He$ in terms of convolution. More precisely, we characterize the group Fourier transform on $L^1(\He,\Hi)$ and $L^2(\He,\Hi)$. In addition, we obtain a similar characterization for the group Fourier transform of $\Hi$-valued Schwartz class functions on the Heisenberg group.
\end{abstract}
\maketitle

\section{Introduction}
The Fourier transform can be characterised as the unique operator satisfying a suitable collection of fundamental properties. This leads naturally to the problem of characterisation: rather than starting with the explicit definition of the Fourier transform, one begins with an operator satisfying certain basic axioms and asks whether these properties uniquely determine the Fourier transform. This is commonly referred to as the axiomatic approach to the characterisation of the Fourier transform.

It is a classical fact that the Fourier transform $\mathcal{F}$ is an isomorphism of $\mathcal{S}(\mathbb{R}^n)$, the Schwartz space of rapidly decreasing functions on $\mathbb{R}^n$, and extends uniquely to an isomorphism of $\mathcal{S}'(\mathbb{R}^n)$, the space of tempered distributions. The problem of characterising classical transforms through minimal structural properties was initiated by the pioneering work of Artstein-Avidan and Milman \cite{artstein2009concept}. They showed, remarkably, that the Legendre transform can be essentially characterised by a few basic properties. This characterisation philosophy was later extended to the classical Fourier transform. In \cite{Alesker_paper}, Alesker, Artstein-Avidan and Milman established that, up to natural normalisations, the Fourier transform is essentially the unique transform on the Schwartz space that interchanges convolution and pointwise multiplication. These results provide a natural motivation for seeking analogous intrinsic characterisations of other Fourier-type transforms. We recall their characterisation theorem below.

\begin{theorem}[\cite{Alesker_paper, Alesker_BC}]\label{th: alesker}
Let $T:S(\mathbb{R}^n)\to S(\mathbb{R}^n)$ be a bijection that admits a bijective extension
\[
T':S'(\mathbb{R}^n)\to S'(\mathbb{R}^n).
\]
Assume that
\begin{align*}
T'(f\ast g) &= T'(f)\cdot T(g),~~ \text{for all } f\in S'(\mathbb{R}^n),\; g\in S(\mathbb{R}^n),\\
T'(f\cdot g) &= T'(f)\ast T(g),~~ \text{for all } f\in S'(\mathbb{R}^n),\; g\in S(\mathbb{R}^n).
\end{align*}
Then there exists $B\in GL_n(\mathbb{R})$ with $|\det(B)|=1$ such that, for every $f\in S(\mathbb{R}^n)$, either
\[
T(f)=\mathcal{F}(f\circ B)
\]
or
\[
T(f)=\overline{\mathcal{F}(f\circ B)}.
\]
\end{theorem}
Later, in \cite{Jaming}, Jaming obtained a characterization of the Fourier transform for $L^1(\R^n)$. Here we state his result.

\begin{theorem}[\cite{Jaming}]\label{th: jaming}
    Let $T:L^1(\R^n)\to C_0(\R^n)$ be a bounded linear operator such that
    \[
    T(f\ast g)=T(f)\cdot T(g),~~\text{ for all } f,g\in L^1(\R^n).
    \]
    Then there exist $E\subset \R^n$ and $\phi:\R^n\to\R^n$ such that
    \[
    T(f)(\xi)=\mathbf{I}_E(\xi)\f(f)(\phi(\xi)), \,\, \xi\in\R^n,
    \]
    holds for all $f\in L^1(\R^n)$, where $\mathbf{I}_E$ is the indicator function of $E$.
\end{theorem}

\par In \cite{Jaming}, the author asked whether such a characterization is possible for the Weyl transform using the the following property.
$$
W_\lambda(f\times_\lambda g)=W_\lambda (f)\circ W_\lambda(g).
$$
The authors in \cite{Lavanya_Weyl},  obtained a characterization of the Weyl transform, which was subsequently employed to establish a characterization theorem for the group Fourier transform on the Heisenberg group in \cite{Lavanya_Heisen}. We now state the characterization theorem for the Weyl transform and for the group Fourier transform on the Heisenberg group. For the definition of convolutions for scalar-valued functions, we refer to \cite{Lavanya_Heisen,Lavanya_Revisiting}.
\begin{theorem}[\cite{Lavanya_Heisen}] \label{T2.6}
    Let $T:\left( L^1(\C^n), \times_\lambda\right)\to \mathcal{B}\left( L^2(\R^n)\right)$ be a nonzero bounded linear map satisfying:
    \begin{itemize}
        \item [(i)] $T(f\times_\lambda g)=T(f)T(g)$, for all $f,g \in L^1(\C^n),$ where $\times_\lambda$ denotes the $\lambda$-twisted convolution (see ),
        \item [(i)] $T(f^\ast)=T(f)^\ast$, where $f^\ast(z)= \overline{f(-z)}.$
    \end{itemize}
    Then there is a subspace $\Hi^\lambda$ of $L^2(\R^n)$ and a unitary representation $\rho_\lambda$ of $\He$ on $\Hi^\lambda$ such that
    $$T(f)= \int_{\C^n} f(z)\rho_\lambda(z,0)\, dz,~ \text{ on } \Hi^\lambda, $$ and there is a decomposition $L^2(\R^n)= \Hi^\lambda\bigoplus V^\lambda$, where 
    $$V^\lambda:=\{v\in L^2(\R^n): (Tf)(v)=0 \text{ for all } f\in L^1(\C^n)\}.$$
\end{theorem}

Appealing to the above theorem, the authors obtained the following
\begin{theorem}[\cite{Lavanya_Heisen}]
    Let $T: L^1(\He)\to L^\infty(\R^\ast, S_2, d\mu)$ be a nonzero continuous linear map satisfying 
    \begin{itemize}
        \item [(i)] $T(f^\ast)(\lambda)=Tf(\lambda)^\ast,~$ for all  $ \lambda\in \R^\ast, f\in L^1(\He),$  where $f^\ast(z,t)= \overline{f(-z,-t)},$
        \item [(ii)] $T(f\ast g)(\lambda)=(Tf)(\lambda)(Tg)(\lambda), ~$ $\lambda\in \R^\ast, ~ f,g \in L^1(\He),$ and
        \item [(iii)] $T(R_{(0,t)}f)(\lambda)= Tf(\lambda)e^{-i\lambda t}, ~ \lambda\in \R^\ast, f\in L^1(\He), t\in \R.$
    \end{itemize}
    Then for each $\lambda\in A$, there is a decomposition $L^2(\R^n)= \Hi^\lambda\bigoplus V^\lambda$, and a unitary representation $\rho_\lambda$ of $\He$ on $\Hi^\lambda$ such that 
    $$Tf(\lambda)= \int_{\He}f(z,t)\rho_\lambda(z,t)\, dz\, dt, \text{ on } \Hi^\lambda,$$
    where $A:=\{\lambda\in \R^\ast: Tf(\lambda)\neq 0 \text{ for some } f\in L^1(\He)\}.$
\end{theorem}

A similar characterization is obtained for the group Fourier transform on $L^2(\He)$. They also obtained a characterization of the Fourier transform form $S(\He)$ onto $\widehat{S}(\He)$, where $S(\He)$ is the Schwartz space on the Heisenberg group and $\widehat{S}(\He)$ is the image of $S(\He)$ under the group Fourier transform. We can view $\widehat{S}(\He)$ as a subspace of $L^2(\R^\ast, S_2, d\mu).$ Here $S_2$ is the space of Hilbert-Schmidt operators on $L^2(\R^n)$ and $L^2(\R^\ast, S_2, d\mu)$ is the space of all functions on $\R^\ast$ taking values in $S_2$, and square integrable with respect to the Plancherel measure $d\mu(\lambda)=c_n|\lambda|^nd\lambda,$ (see \cite{geller1977fourier}). Now we recall their result below.

\begin{theorem}[\cite{Lavanya_Revisiting}]
    Let $T: S(\He)\to \widehat{S}(\He)$ be a bijection which satisfies 
    \begin{itemize}
        \item [(i)] $T(f\ast g)= (Tf) (Tg)$, and
        \item [(ii)] $T(f\ast_3 g)(\lambda)= 2\pi \int_{\C^n} \f^{-1} f(w, \lambda)Tg(\lambda, w)\, dw.$
    \end{itemize}
    Then there exists a map $\zeta: \R^\ast\to \C^n$ such that $Tf(\lambda)= \widehat{f}(\lambda, \zeta(\lambda)).$
\end{theorem} 

\par Further, the characterization of the Fourier transform on compact groups and on locally compact abelian groups is investigated in \cite{Kumar} and \cite{Krukowski_RIMA} respectively. Furthermore, in \cite{krukowski}, the author has demonstrated that two fundamental properties of the Fourier transform, namely the time differentiation property and the Dirac delta property, are sufficient to characterize the Fourier transform among integral transforms.

\par 
On the other hand, in \cite{Hasan1}, the characterization problem is investigated for the Fourier transform of Hilbert space-valued functions on $\R$. Let $\Hi$ be a separable Hilbert space. Then, for
$f\in L^1(\R,\Hi),$ the Fourier transform $\widehat{f}$ of a vector valued function $f$ is defined as in the scalar valued case:
\begin{equation*}
    \widehat{f}(\xi)=\int_\R e^{-2\pi ix\xi}f(x)dx,~~\xi\in \R,
\end{equation*}
where the right hand side is a Bochner integral. For a detailed study of integration of vector-valued functions, we refer to \cite{Diestel}. 

Let us now state the characterization theorem in the context of vector valued functions. For the definition of two convolutions for Hilbert space valued functions, we refer to \cite{Hasan1}. See also \cite{AHM}.

\begin{theorem}[\cite{Hasan1}]\label{L1 char}
Let $(\Hi, \langle {\cdot, \cdot } \rangle)$ be a separable Hilbert space with an orthonormal basis $\{e_k: k\in\N\}$ and $T:L^1(\R, \Hi)\rightarrow C_0(\R, \Hi)$ be a bounded linear operator satisfying
\item[(i)] $T((f\ast g)\odot h)=\left(Tf\times Tg\right) \smblksquare Th$, for $f, g, h \in L^1(\R, \Hi)$,
\item[(ii)]$T(\tau_s f)(\xi)=e^{-2\pi i s \xi}Tf(\xi)$, $f\in L^1(\R, \Hi), ~~s, \xi \in \R$,
\item[(iii)] If $\langle {f(x), e_k}\rangle =0$ for $a.e. ~x \in \R$ for some $k$ then $\langle {Tf(\xi), e_k}\rangle =0$ for all $ \xi \in \R$, and
\item[(iv)] $T(\widetilde {f}~)(-\xi)=Tf(\xi)$,  $f\in L^1(\R, \Hi),~~ \xi \in \R$, where $\widetilde {f}(x)=f(-x)$.\\
Then there exists $E\subset \R$ such that
\begin{center}
$\langle {Tf(\xi), Tg(\xi)}\rangle =\mathbf{I}_E(\xi)\langle {\widehat{f}(\xi), \widehat{g}(\xi)}\rangle$,\quad  $f,g \in L^1(\R, \Hi), ~ \xi \in \R$,
\end{center}
where $\mathbf{I}_E$ is the indicator function of $E$.
\end{theorem}
As a corollary to Theorem \ref{L1 char}, the following characterization is obtained.
\begin{corollary}[\cite{Hasan1}]\label{cor: L^1 cor}
If in addition to the hypothesis of Theorem \ref{L1 char}, $T$ satisfies $T(g_k)=g_k$ where $g_k(x)=e^{-\pi x^2}e_k$ for all $k\in \N$, then there exists $E\subset \R$ such that $Tf(\xi)=\chi_E (\xi) \widehat{f}(\xi)$,\quad $f \in L^1(\R, \Hi),~~\xi \in \R$.
\end{corollary}

In this work, we establish such a characterization for the Fourier transform on the Heisenberg group $\He$ within the Hilbert space framework. Let $\Hi$ be a Hilbert space and $f\in L^1(\He, \Hi)$ then the group Fourier transform of $f$ is defined as (see \cite{Swain})
$$\widehat{f}(\lambda)(\phi)(\xi)=\int_{\He}f\left(z,t\right)\pi_\lambda(z,t)\phi(\xi)\, dz\, dt, ~~~~ \lambda\in \R^*,\,\phi\in L^2(\R^n).$$

More precisely, we consider the Fourier transform on $L^1(\He,\Hi)$, $L^2(\He,\Hi)$, $S(\He,\Hi)$ and obtain respective characterization. We first define the notion of convolution for Hilbert space valued functions on the Heisenberg group as follows.
$$f\ast g(z,t)= \int_{\He}\left\langle f(w,s), g\left(\left(z,t\right)\left(-w,-s\right)\right)\right\rangle\, dw\, ds, ~~~~ (z,t)\in \He, 
$$
where $\langle\cdot,\cdot\rangle$ is the inner product on $\Hi$. As convolution of two vector valued functions is a scalar valued function, we defined convolution of three functions and obtain the corresponding convolution theorem (see \eqref{2 conv}). 
\par Following \cite{Lavanya_Heisen}, we make use of the characterization of the Weyl transform to obtain the characterization of the Fourier transform on $L^1(\He,\Hi)$. We use Theorem \ref{T2.6} instead of proving its vector valued analogue. In order to use the characterization theorem for the Weyl transform of scalar valued functions, we use Cohen's factorization theorem for Banach modules. 
\par By viewing $L^2(\He, \Hi)$ as the tensor product of $L^2(\He)$ and $\Hi$, we observe that the Fourier transform on $L^2(\He,\Hi)$ has a special form. More precisely, for $f\in L^2(\R^n),\,u\in\Hi$, the Fourier transform of the basic tensor $(f\otimes u)$ is
$$
(f\otimes u)~\widehat{}~(\lambda)=(\widehat{f}(\lambda)\circledast u),~\text{ where } (T\circledast u)\phi(\xi)= T\phi(\xi)\,u,~\,T\in \mathcal{B}\left(L^2\left(\R^n\right)\right)
$$
So, we obtain characterization for operators defined on $L^2(\He,\Hi)$ having the same form as the Fourier transform. We prove a multiplier type theorem on $L^2(\He, \Hi)$ and use it to get the desired characterization of the Fourier transform on $L^2(\He,\Hi)$.
\par In \cite{Lavanya_Revisiting}, the authors considered the following unitary representations, parametrized by $(\lambda,\zeta)\in\R^\ast\times\C^n$, 
$$
\pi_{\lambda,\zeta}(z,t)=\chi_\zeta(z,t)\pi_\lambda(z,t)
$$
where $\chi_\zeta(z,t)=e^{i\Re(z\cdot\bar{\zeta})}$, are one dimensional representations of the Heisenberg group and defined the Fourier-Weyl transform as
$$
\F(f)(\lambda,\zeta):=\int_{\He}f(z,t)\pi_{\lambda,\zeta}(z,t)dzdt.
$$
The authors in \cite{Lavanya_Revisiting} obtained a characterization for the Fourier-Weyl transform and subsequently, characterized the Fourier transform on the Heisenberg group.
For $f\in L^1(\He,\Hi)$, we define the Fourier-Weyl transform of vector valued functions as
$$
(\F f(\lambda,\zeta)\phi)(\xi):=\int_{\He}f(z,t)\pi_{\lambda,\zeta}(z,t)\phi(\xi) dzdt,
$$
$\phi\in L^2(\R^n)$ and establish that the behavior of the Fourier-Weyl transform on this convolution and the usual convolution characterizes it. As a consequence of this, we obtain the characterization of the Fourier transform on $S(\He,\Hi)$.

\section{Preliminaries} \label{prelim}

We begin by recalling that the Heisenberg group $\He$ is  $\C^n\times \R$ equipped with the group law 
$$(z,t)\cdot(w,s)=\left(z+w, t+s+\frac{1}{2}\Im\left(z.\overline{w}\right)\right)$$
The Haar measure on $\He$ is $dz\,dt$, the Lebesgue measure on $\C^n\times\R.$ Stone-von Neumann theorem provides all the infinite-dimensional irreducible unitary representations of $\He$. These are parametrized by $\lambda\in \R^*$ and their action on $L^2(\R^n)$ is given by 
$$\pi_\lambda(z,t)\varphi(\xi)= e^{i\lambda t}e^{i\lambda\left(x.\xi+ \frac{1}{2}x.y\right)}\varphi\left(\xi+y\right),~~~~~~\xi\in \R^n, ~~~~\varphi\in L^2(\R^n).$$
The group Fourier transform of $f\in L^1(\He)$ is defined by
$$
\widehat{f}(\lambda)=\int_{\He} f(z,t)\pi_\lambda(z,t)\,dzdt.
$$
Let
$$
f^\lambda(z)=\int_{-\infty}^\infty f(z,t)e^{i\lambda t}dt
$$
be the inverse Fourier transform of $f$ in the central variable, then $\widehat{f}(\lambda)=W_\lambda(f^\lambda)$, where $W_\lambda$, known as the Weyl transform, is defined as
$$
W_\lambda(f)=\int_{\C^n} f(z)\pi_\lambda(z) dz,~~\,\pi_\lambda(z):=\pi_\lambda(z,0).
$$
We define the convolution of $f,g\in L^1(\He)$ in terms of the left translation as
$$
(f\ast g)(z,t)=\int_{\He}f(w,s)\,g((-w,-s)\cdot(z,t))\,dwds.
$$
The convolution for functions on the Heisenberg group induces a convolution for functions on $\C^n$ as follows.
$$
(f\ast g)^\lambda= (f^\lambda\times_\lambda g^\lambda),
$$
where
\[F\times_\lambda  G(z)= \int_{\C^n} F(z-w) G(w) e^{i\frac{\lambda}{2}\Im\left(z\cdot \overline{w}\right)}\, dw, ~~~~ z\in \C^n,\]
is called the $\lambda$-twisted convolution. We refer to \cite{Heisen} for a study of harmonic analysis on the Heisenberg group.

\par
In the remaining part of this section we extend the notion of convolution (resp. twisted convolution) of vector valued functions on the Heisenberg group (resp. on $\R^{2n}$) and obtain the convolution theorem. Further we state factorization results which will be useful in the sequel.\\

\par Let $\mathcal{B}\left(L^2(\R^n), L^2(\R^n,\Hi)\right)$ denote the space of all bounded linear maps from $L^2(\R^n)$ to $L^2(\R^n,\Hi)$. Then we can show that $\widehat{f}(\lambda)\in \mathcal{B}\left(L^2(\R^n), L^2(\R^n,\Hi)\right),$ with $\|\widehat{f}(\lambda)\|_{op}\leq \|f\|_{L^1(\He,\Hi)}.$

Recall that the convolution $f,g \in L^1(\He,\Hi)$ is defined by
$$f\ast g(z,t)= \int_{\He}\left\langle f(w,s), g\left(\left(z,t\right)\left(-w,-s\right)\right)\right\rangle\, dw\, ds, ~~~~ (z,t)\in \He, 
$$
where $\langle\cdot,\cdot\rangle$ is the inner product on $\Hi$, which is linear in the first variable and conjugate linear in the second variable. Observe that the convolution of two vector valued functions on the Heisenberg group turns out to be a scalar valued function. So,
for $\varphi\in L^1(\He)$ and $h\in L^1(\He,\Hi)$, we define
\begin{equation}\label{2 conv}
\varphi\fisheye h(z,t)=\int_{\He}\varphi(w,s)h\left(\left(z,t\right)\left(-w,-s\right)\right)\, dw\, ds, ~~~~~(z,t)\in \He. 
\end{equation}
Then it is easy to see that $(f\ast g)\fisheye h\in L^1(\He,\Hi)$, whenever $f,g,h\in L^1(\He,\Hi)$. These convolutions on $\He$ induces $\lambda$-twisted convolution on $\R^{2n}$. In fact, if
$$
f^\lambda(z)=\int_{-\infty}^\infty f(z,t)e^{i\lambda t}\, dt, ~~~~~ z\in \C^n
$$
is the inverse Fourier transform of $f\in L^1(\He,\Hi)$ in the central variable then it is easy to see that 
$$\left( f_1\ast f_2\right)^\lambda= f_1^\lambda \ast_\lambda f_2^{-\lambda},$$
where
$$g_1\ast_\lambda g_2(z)= \int_{\C^n}\left\langle g_1(w), g_2(z-w)\right\rangle e^{i\frac{\lambda}{2}\Im\left(z\cdot \overline{w}\right)}\, dw, ~~~~ z\in \C^n.$$
Moreover, 
$$\left( f_1\fisheye f_2\right)^\lambda= f_1^\lambda \fisheye_\lambda f_2^{\lambda},$$
where
$$g_1\fisheye_\lambda g_2(z)=\int_{\C^n}g_1(w)g_2(z-w)e^{i\frac{\lambda}{2}\Im\left(z\overline{w}\right)}\, dw, ~~~~ z\in \C^n.$$

We now obtain the convolution theorem. Observe that, for $g\in L^1(\He,\Hi)$, $\widehat{g}(\lambda)^\ast\in \mathcal{B}\left( L^2(\R^n,\Hi), L^2(\R^n)\right)$. Further, a simple calculation shows that
$$\widehat{g}(\lambda)^\ast \phi(\xi)= \int_{\He} \left\langle \phi(\xi-y), g(z,t)\right\rangle e^{-i\lambda t}e^{-i\lambda \left(x.\xi-\frac{1}{2}x.y\right)}\, dz\, dt \quad \text{ for all } \phi\in L^2(\R^n,H),\, \xi\in \R^n.$$
Therefore,
$$
\widehat{g}(\lambda)^\ast\circ \widehat{f}(\lambda)\phi(\xi)
    = \int_{\He} \left\langle \widehat{f}(\lambda)\phi(\xi-y), g(z,t)\right\rangle e^{-i\lambda t}e^{-i\lambda \left(x.\xi-\frac{1}{2}x.y\right)}\, dz\, dt.
$$
It is easy to verify that the right hand side is the same as
$$
\int_{\He}\int_{\He} \left\langle f(w,s), g(-z,-t)\right\rangle     \pi_\lambda(z,t)\pi_\lambda(w,s)\phi(\xi)\,dw\,ds\,dz\,dt.
$$
On the other hand,
$$
\left(f\ast g\right)~\widehat{}~(\lambda)(\phi)(\xi)
    = \int_{\He}\int_{\He} \left\langle f(w,s), g(z,t)\right\rangle \pi_\lambda(z,t)\pi_\lambda(w,s)\phi(\xi)\, \,dw\,dsdz\, dt
$$
Thus
$$\left(f\ast g\right)~\widehat{}~ (\lambda)= \left(\widetilde{g}\right)~\widehat{}~(\lambda)^\ast\circ \widehat{f}(\lambda), \quad \lambda\in \R^\ast.$$
A similar calculations shows that
$$\left(\varphi \fisheye h\right)~\widehat{}~(\lambda)= \widehat{h}(\lambda)\circ \widehat{\varphi}(\lambda) \quad \text{ for all } \lambda\in \R^\ast.$$
Combining the above equalities we get
$$\left(\left(f\ast g\right)\fisheye h\right)~\widehat{}~(\lambda)= \widehat{h}(\lambda)\circ \left(\left(\widetilde{g}\right)~\widehat{}~(\lambda)^\ast\circ \widehat{f}(\lambda)\right) \text{ for all } \lambda\in \R^\ast,$$
for $f,g,h\in L^1({\He},\Hi)$.\\
We shall characterize the Fourier transform in terms of the above convolution theorem. We also need the following equality of convolutions.

\begin{lemma}\label{lem: conv rel}
    Let $g_1,g_2, g_3, g_4\in L^1(\C^n,\mathcal{H})$. Then $\left( g_1\ast_\lambda g_3\right)\times_\lambda \left(g_2\ast_\lambda g_4\right)=\left(\left(g_2\ast_\lambda g_4\right)\fisheye_\lambda g_1\right)\ast_\lambda g_3$.
\end{lemma}

\begin{proof}
We have 
$$
\left( g_1\ast_\lambda g_3\right)\times_\lambda \left(g_2\ast_\lambda g_4\right)\left(z\right)
        =\int_{\C^n} \left( g_1\ast_\lambda g_3\right)(z-w)\left( g_2\ast_\lambda g_4\right)\left( w\right)e^{i\frac{\lambda}{2}\im \left(z\overline{w}\right)}\, dw.
$$
But the right hand side simplifies to
\begin{align*}
&\int_{\C^n}\int_{\C^n}\int_{\C^n} \left\langle g_1\left(w'\right), g_3\left(z-w-w'\right)\right\rangle \left\langle g_2\left(w''\right), g_4\left(w-w''\right)\right\rangle\\
&\times e^{i\frac{\lambda}{2}\im \left(\left(z-w\right)\overline{w'}\right)}e^{i\frac{\lambda}{2}\im \left(w\overline{w''}\right)}e^{i\frac{\lambda}{2}\im \left(z\overline{w}\right)}\,dw\,dw'\,dw''.
\end{align*}
  
On the other hand, 
\begin{align*}
&\left(\left(g_2\ast_\lambda g_4\right)\fisheye_\lambda g_1\right)\ast_\lambda g_3\left( z\right)= \int_{\C^n} \left\langle \left( g_2\ast_\lambda g_4\right)\fisheye_\lambda g_1\left( w\right), g_3\left( z-w\right)\right\rangle e^{i\frac{\lambda}{2}\im\left(z\overline{w}\right)}\, dw.
\end{align*}
  The right hand side simplifies to
\begin{align*}
  &\int_{\C^n}\int_{\C^n}\int_{\C^n} \left\langle g_2\left( w''\right), g_4\left( w'-w''\right)\right\rangle \left\langle g_1\left( w-w'\right), g_3\left( z-w\right)\right\rangle\\ & \times e^{i\frac{\lambda}{2}\left\{\im \left( w'\overline{w''}\right)+ \im \left( w\overline{w'}\right)+\im \left( z\overline{w}\right)\right\}}\, dw\, dw'\, dw'',
  \end{align*}
  which is same as
  \begin{align*}
      &\int_{\C^n}\int_{\C^n}\int_{\C^n} \left\langle g_2\left( w''\right), g_4\left( w-w''\right)\right\rangle \left\langle g_1\left( w'\right), g_3\left( z-w-w'\right)\right\rangle\\
        &e^{i\frac{\lambda}{2}\{\im \left( w\overline{w''}\right)+ \im \left( (w+w')\overline{w}\right)+\im \left( z\overline{(w+w')}\right)\}}\, dw'\, dw\, dw''
  \end{align*}
  Since, $\im \left( w\overline{w''}\right)+ \im \left( (w+w')\overline{w}\right)+\im \left( z\overline{(w+w')}\right)=\im \left( w\overline{w''}\right)+ \im \left( z\overline{w}\right)+\im \left( (z-w)\overline{w'}\right)$, the desired equality follows.
\end{proof}
Now we state the following factorization theorem, which we will use to prove the characterization theorem for $L^1(\He,\Hi).$ 
\begin{theorem}[Cohen Factorization Theorem]
    If $A$ is a Banach algebra with a bounded left approximate identity, then every element $g\in A$ can be factored as $g=g_1.g_2$ for some $g_1,g_2\in A.$
\end{theorem}
Let $g\in L^1(\C^n)$. Then applying the above theorem to the Banach algebra $(L^1(\C^n),\ast_\lambda)$ we get $g=g_1\ast_\lambda g_2$, for some $g_1,g_2\in L^1(\C^n).$ Let $G_1=g_1\otimes e_1$ and $G_2=\bar{g_2}\otimes e_1$, then 
\begin{align*}
    G_1\ast_\lambda G_2(z) = \int_{\C^n} \left\langle G_1(w),G_2(z-w)\right\rangle e^{i\frac{\lambda}{2}\Im(z\overline{w})}\,dw
    = g_1\ast_\lambda g_2(z)
    = g(z).
\end{align*}
Therefore any element $g\in L^1(\C^n)$, can be written as $g=G_1\ast_\lambda G_2$ for some $G_1,G_2\in L^1(\C^n,\Hi).$ More generally, for $f, g\in L^1(\C^n)$, then there exist $F,G,H \in L^1(\C^n,H)$ such that $f= H\ast_\lambda F$ and $g= H\ast_\lambda G$. This follows from the following factorization theorem.

\begin{theorem}\label{th: CW fact}[Cohen- Hewitt Factorization Theorem]
    Let $A$ be a Banach algebra with a bounded left approximation identity $\{e_\alpha\}_{\alpha\in I}$ and $X$ be a left Banach $A$-module. If $x\in X$ can be approximated by elements of the form $e_\alpha.x$, i.e. $\lim_\alpha e_\alpha.x=x$, then there exist $y\in X$ and $a\in A$ such that $x=a.y$.
\end{theorem}
Let $A$ denote the Banach algebra $\left(L^1(\C^n), \ast_\lambda\right)
$. Consider $X=A\times A$ equipped with the norm $$\|(f,g)\|_X:=\max \left\{\|f\|_A,\|g\|_A\right\}$$ 
and define 
$$h.(f,g):=(h\ast_\lambda f,h\ast_\lambda g),~~ \text{ for }h\in A~\text{ and }(f,g)\in X.$$
Then $X$ is a Banach $A$-module. Moreover, if $\{e_\alpha\}$ is a bounded left approximate identity in $A$, then it is easy to check that $\lim_\alpha e_\alpha.(f,g)=(f,g)$ for all $(f,g)\in X.$\\
Let $(f,g)\in X$. Then, by making use of the Theorem \ref{th: CW fact}, we have $h\in A$ and $(f',g')\in X$ such that 
$$(f,g)=(h\ast_\lambda f', h\ast_\lambda g').$$
Now, setting $F=f'\otimes e_1, G=g'\otimes e_1, H= h\otimes e_1$ , it follows that $f= H\ast_\lambda F$ and $g=H\ast_\lambda G$. Note that here the last two convolutions are twisted convolution for Hilbert space valued functions.

We state this observation in the following

\begin{corollary} \label{c2.5}
    Let $f, g\in L^1(\C^n)$, then there exists $F,G,H \in L^1(\C^n,\Hi)$ such that $f= H\ast_\lambda F$ and $g= H\ast_\lambda G$.
\end{corollary}

\section{Characterization on $L^1(\He,\Hi)$ }
In this section we prove a characterization for the group Fourier transform of vector valued functions on the Heisenberg group. Towards this end, appealing to the characterization of the Weyl transform of scalar valued functions, we derive the following result. Subsequently, we obtain the desired characterization. We will use the following lemma in the proof.
\begin{lemma} \label{L2.2}
    Let $A$ be a Banach algebra with a bounded left approximate identity. Let $\{z_n\}$ be a sequence in $A$ such that $\lim_n z_n=0$. Then there exist $a\in A$ and a sequence $\{y_n\}$ in $A$ such that $z_n=a.y_n$ for all $n$ and $\lim_n y_n=0.$
\end{lemma}

\begin{theorem}
    Let $T: L^1(\He, \Hi)\to L^\infty\left(\R^\ast, \mathcal{B}\left(L^2\left(\R^n\right), L^2\left(\R^n, \Hi\right)\right),d\mu\right)$ be a nonzero continuous linear map satisfying:
    \begin{itemize}
        \item [(i)] $T\left( \left( f\ast g\right)\fisheye h\right)\left(\lambda\right)= Th(\lambda)\circ  T\widetilde{g}\left(\lambda\right)^\ast \circ Tf\left(\lambda\right),~~~~~~~$ for all $f,g,h \in L^1(\He, \Hi ).$
        \item [(ii)] $T\left( R_{(0,t)}f\right)\left( \lambda\right)= e^{-i\lambda t}Tf\left(\lambda \right) $, $~~~~~\lambda \in \R^\ast$, $f\in L^1(\He, \Hi )$, $t\in \R.$
    \end{itemize}
    Then for each $\lambda \in A$, there is a decomposition of $L^2(\R^n)= \mathcal{H}^\lambda\oplus V^\lambda$, and a unitary representation $\rho_\lambda$ of $\He$ on $\mathcal{H}^\lambda$ such that, 
    $$T(\Tilde{g})(\lambda)^\ast \circ Tf(\lambda)=\int_{\He}(f\ast g)(\lambda)\rho_\lambda (z,t)\, dz\, dt, ~~~~~~ \text{ on } \mathcal{H}^\lambda,$$ where $A:= \{\lambda \in \R^\ast: Tf(\lambda)\neq 0 ~~\text{ for some } f\in L^1(\He, \Hi)\}.$
\end{theorem}
\begin{proof}
    For $\lambda \in \R^\ast$, let $T_\lambda (f):= Tf(\lambda)$, for all $f\in L^1(\He,\mathcal{H}).$ 
    For fixed $\phi\in L^2(\R^n)$ and $\psi \in L^2(\R^n,\Hi),$ we consider the map $f\mapsto \left( T_\lambda\left(f\right)\phi, \psi\right)$, where $$\left( T_\lambda\left(f\right)\phi, \psi\right)= \int_{\R^n}\left\langle T_\lambda (f)\phi(x), \psi(x)\right\rangle \, dx.$$ Then, appealing to the Cauchy-Schwartz inequality for $L^2(\R^n,\Hi)$ and the boundedness of $T$, it follows that the map $f\mapsto \left( T_\lambda\left(f\right)\phi, \psi\right)$ is a linear continuous functional on $L^1(\He,\Hi).$ In fact, 
    \begin{align*}
        &|(T_\lambda(f)\phi, \psi)| \leq \| T_\lambda (f)\phi\|_{L^2(\R^n,H)}\|\psi\|_{L^2(\R^n,H)}\leq \|T_\lambda(f)\|\|\phi\|_{L^2(\R^n)}\|\psi\|_{L^2(\R^n,H)}\\
        &\leq C\|f\|_{L^1(\He,H)}\|\phi\|_{L^2(\R^n)}\|\psi\|_{L^2(\R^n,H)}.
    \end{align*}
    Then there exists a $F_\lambda \in L^\infty(\He, \Hi )$ such that 
    $$\left(T_\lambda\left(f\right)\phi, \psi\right)= \int_{\He}\left\langle f(z,t), F_\lambda\left(\left(z,t\right);\phi, \psi\right)\right\rangle \, dz\, dt, ~~~~~~ ~~~~ f\in L^1(\He,\Hi). $$
    Now, let $f\in L^1(\He,\Hi)$ be of the form $f(z,t)= g(z)h(t)$ where $g\in L^1(\C^n,\Hi)$ and $h\in L^1(\R).$ Then
    \begin{align*}
    \left( T_\lambda\left(f\right)\phi, \psi\right)&= \int_{\He}\left\langle f(z,t), F_\lambda\left(\left( z,t\right);\phi, \psi\right)\right\rangle \, dz\,dt\\
    &=\int_{-\infty}^\infty h(t)\left(\int_{\C^n}\left\langle g(z), F_\lambda\left(\left( z,t\right);\phi, \psi\right)\right\rangle \, dz\right)\,dt\\
    &= \int_{-\infty}^\infty h(t)\Phi_\lambda(t)\,dt,
    \end{align*}
    where $\Phi_\lambda(t)=\int_{\C^n}\left\langle g(z), F_\lambda\left(\left( z,t\right);\phi,\psi\right)\right\rangle\,dz.$ Moreover, by making use of the hypothesis (ii), for any $s\in\R$, we get
    $$\left( T_\lambda\left(f\right)e^{-i\lambda s}\phi, \psi\right)= \left( T_\lambda \left(R_{\left(0,s\right)}f\right)\phi,\psi\right)
    =\int_{-\infty}^\infty h(t)\Phi_\lambda(t-s)\,dt.$$
    Thus $\Phi_\lambda(t-s)= e^{-i\lambda s}\Phi_\lambda(t)$ for all $s\in \R$, $a.e. \,\,t\in \R.$ Let $\Psi$ be a Schwartz class function satisfying $\widehat{\Psi}(\lambda)\neq 0.$ Then
    $$\int_{-\infty}^\infty\Phi_\lambda(t-s)\Psi(s)\,ds= \int_{-\infty}^\infty e^{-i\lambda s}\Phi_\lambda(t)\Psi(s)\,ds=\Phi_\lambda(t)\widehat{\Psi}(\lambda).$$
    As the left hand side is a smooth function of t, so is $\Phi_\lambda.$ Therefore we have, $\Phi_\lambda(t-s)= e^{-i\lambda s}\Phi_\lambda(t)$ for all $s,t\in \R$, which implies that $\Phi_\lambda(t)= e^{i\lambda t}\Phi_\lambda(0)$ for all $t\in \R.$ Thus for every $g\in L^1(\C^n,\Hi),$ 
    $$t\mapsto \int_{\C^n}\left\langle g(z), F_\lambda\left(z,t\right);\phi,\psi\right\rangle\, dz$$ 
    is a continuous function and
    $$\int_{\C^n}\left\langle g(z), F_\lambda\left(\left(z,t\right);\phi,\psi\right)\right\rangle\, dz=\Phi_\lambda(t)=e^{i\lambda t}\Phi_\lambda(0)= e^{i\lambda t}\int_{\C^n}\left\langle g(z), F_\lambda\left(\left(z,0\right);\phi,\psi\right)\right\rangle\, dz.$$
    {\color{black} Let $\{e_k:k\in\N\}$ be an orthonormal basis for $\Hi$. Then, by choosing $g_k(z)=G(z)e_k$, for some $G\in L^1(\C^n)$, it follows that
    $$
    \int_{\C^n}G(z)\left\langle e_k, F_\lambda\left(\left(z,t\right);\phi,\psi\right)\right\rangle\, dz= e^{i\lambda t}\int_{\C^n}G(z)\left\langle e_k, F_\lambda\left(\left(z,0\right);\phi,\psi\right)\right\rangle\, dz,
    $$
    for all $G\in L^1(\C^n)$. Hence, using the $L^1-L^\infty$ duality, it follows that
    $$\left\langle e_k, F_\lambda\left(\left(w,t\right);\phi,\psi\right)\right\rangle= e^{i\lambda t}\left\langle e_k, F_\lambda\left(\left(w,0\right);\phi,\psi\right)\right\rangle.$$
    Since $k$ is arbitrary
    $$F_\lambda\left(\left(w,t\right);\phi,\psi\right)= e^{-i\lambda t}F_\lambda\left(\left(w,0\right);\phi,\psi\right), ~~~~~~~~ a.e.~~~ w\in \C^n.$$
    }
    Then
    \begin{align*}
        \left( T_\lambda\left(f\right)\phi,\psi\right) &= \int_{\He}\left\langle f(z,t), e^{-i\lambda t} F_\lambda\left(\left(z,0\right);\phi,\psi\right)\right\rangle\, dz\,dt\\
        &=\int_{\He}\left\langle e^{i\lambda t}f(z,t),  F_\lambda\left(\left(z,0\right);\phi,\psi\right)\right\rangle\, dz\,dt\\
        &=\int_{\C^n}\left\langle \int_{-\infty}^\infty f(z,t)e^{i\lambda t}\,dt,  F_\lambda\left(\left(z,0\right);\phi,\psi\right)\right\rangle\, dz\\
        &=\int_{\C^n}\left\langle f^\lambda(z),  F_\lambda\left(\left(z,0\right);\phi,\psi\right)\right\rangle\, dz,
    \end{align*}
    which shows that $T_\lambda(f)$ depends only on $f^\lambda$ and satisfies
    $$\|T_\lambda(f)\|\leq C\|f^\lambda\|_{L^1(\C^n,\Hi)}.$$
    Now for given $\lambda\in \R^*$, fix $\psi \in L^1(\R)$ such that $\widehat {\psi}(-\lambda)=1$. Then we define an operator {\color{black}$S_\lambda: L^1(\C^n)\to 
\mathcal{B}(L^2(\R^n))$} by
    $$S_\lambda(g):= \left(T_\lambda\left(\widetilde{g_2}\psi\right)\right)^\ast\circ T_\lambda\left(g_1\psi\right),$$
    where $g=g_1\ast_\lambda g_2,$ for some $g_1,g_2\in L^1(\C^n,\mathcal{H})$, this is possible using the factorization theorem. In order to see the well-definedness of $S_\lambda$,
let 
$ g_1\ast_\lambda  g_2= g_3\ast_\lambda g_4$. Then $\left(\left( g_1\ast_\lambda g_2\right)\fisheye_\lambda h\right)\psi= \left(\left( g_3\ast_\lambda g_4\right)\fisheye_\lambda h\right)\psi,$ for all $h\in L^1(\C^n)$. Since,
    \begin{align*}
    \left(\left( g_1\psi\ast g_2\widetilde{\psi}\right)\fisheye h\psi\right)^\lambda(z)
        =\left( g_1\ast_\lambda g_2\right)\fisheye_\lambda h(z)= \left(\left(( g_1\ast_\lambda g_2\right)\fisheye_\lambda h)\psi\right)^\lambda(z),
    \end{align*}
    and $T_\lambda(f)$ depends only on $f^\lambda$, it follows that
    $$
    T_\lambda\left(\left( g_1\psi\ast g_2\widetilde{\psi}\right)\fisheye h\psi\right)=T_\lambda\left(\left( g_3\psi\ast g_4\widetilde{\psi}\right)\fisheye h\psi\right)
    $$
    for all $ h\in L^1(\C^n)$. Now, appealing to the hypothesis (i), it is easy to see that $S_\lambda( g_1\ast_\lambda g_2)=S_\lambda( g_3\ast_\lambda g_4)$.

We now show that $S_\lambda$ is a $\ast$-homomorphism from $L^1(\C^n)$ into 
{\color{black}$\mathcal{B}(L^2(\R^n))$}.\\
It is easy to see that $S_\lambda$ is linear. In fact, if $f,g\in L^1(\C^n)$, then by making use of the Corollary \ref{c2.5}, choose $F,G,H\in L^1(\C^n,H)$ such that $f=H\ast_\lambda F$ and $g=H\ast_\lambda G$. Then from the definition it follow that $S_\lambda(\alpha f+g)=S_\lambda(f)+S_\lambda(g),$ $\alpha\in\C$.\\
In order to see that $S_\lambda$ is continuous, 
take a sequence $\{g_n\}$ in $L^1(\C^n)$ such that $g_n\to 0$. Then appealing to Lemma \ref{L2.2}, choose $h,f_n\in L^1(\C^n)$ such that $g_n= h\ast_\lambda f_n$ and $\lim_n f_n=0.$ Set $H=h\otimes e_1$ and $F_n=f_n\otimes e_1$. Then, using the continuity of $T$, 
$$
S_\lambda(g_n)= S_\lambda \left( H\ast_\lambda F_n\right)= \left(T_\lambda\left( \widetilde{F_n}\psi\right)\right)^\ast\circ T_\lambda\left(H\psi\right)\to 0,
$$
which shows that $S_\lambda$ is continuous.

We now show that $S_\lambda(G_1\times_\lambda G_2)= S_\lambda(G_1)\circ S_\lambda(G_2),$ for all $G_1,\,G_2\in L^1(\C^n)$. Choose $g, g_i\in L^1(\C^n,\Hi)$ such that $G_1=  g\ast g_1$ and $G_2= g\ast g_2$. Then the definition of $S_\lambda$ gives
$$
S_\lambda(G_1)\circ S_\lambda(G_2)= S_\lambda( g\ast g_1)\circ S_\lambda( g\ast g_2)
        = \left(T_\lambda\left( \widetilde{ g_1}\psi\right)^\ast \circ T_\lambda\left( g\psi\right)\right)\circ \left(\left(T_\lambda\left( \widetilde{ g_2}\psi\right)\right)^\ast \circ T_\lambda\left( g\psi\right)\right).
$$
Now using the hypothesis (i), $$T_\lambda\left( g\psi\right)\circ \left(T_\lambda\left( \widetilde{ g_2}\psi\right)\right)^\ast \circ T_\lambda\left( g\psi\right)=T_\lambda\left(\left( g\psi\ast  g_2\widetilde{\psi}\right)\fisheye  g\psi\right).$$ 
But $\left( g\psi\ast  g_2\widetilde{\psi}\right)\fisheye  g\psi=\left(\left(g\ast_\lambda  g_2\right)\fisheye_\lambda g\right)\psi$. Thus
$$
S_\lambda(G_1)\circ S_\lambda(G_2)= S_\lambda\left(\left(\left( g\ast_\lambda  g_2\right)\fisheye_\lambda g\right)\ast_\lambda g_1\right)= S_\lambda\left(\left( g\ast_\lambda g_1\right)\times_\lambda \left( g\ast_\lambda g_2\right)\right)= S_\lambda\left(G_1\times_\lambda G_2\right),
$$
using Lemma \ref{lem: conv rel}. Finally,
    \begin{align*}
        \left(S_\lambda\left( g_1\ast_\lambda g_2\right)\right)^\ast = \left(T_\lambda\left( g_1\psi\right)\right)^\ast\circ T_\lambda\left(\widetilde{g_2}\psi\right)=S_\lambda\left( \widetilde{g_2}\ast_\lambda \widetilde{g_1}\right)=S_\lambda\left(\left( g_1\ast_\lambda g_2\right)^\ast\right),
    \end{align*}
    shows that $S_\lambda$ is a $\ast$-homomorphism. Thus, appealing to Theorem \ref{T2.6}, there exists a representation $\rho_\lambda$ of $\He$ on $\Hi^\lambda$ such that
    \begin{align*}
        \left.S_\lambda\left( \left( f_1\ast f_2 \right)^\lambda\right)\right|_{\mathcal{H}^\lambda}&= \int_{\C^n} \left(f_1\ast f_2\right)^\lambda(z)\rho_\lambda(z,0)\,dz\\
        &=\int_{\C^n}\int_{-\infty}^\infty\left(f_1\ast f_2\right)\left(z,t\right)e^{i\lambda t}\rho_\lambda(z,0)\,dt\,dz\\
        &=\int_{\He}\left(f_1\ast f_2\right)\left(z,t\right)\rho_\lambda\left(z,t\right)\,dz\,dt,~~\text{ on }~\mathcal{H}^\lambda,
    \end{align*}
   where $L^2(\R^n)=\Hi^\lambda\bigoplus V^\lambda$, and 
    $$V^\lambda:=\left\{v\in L^2(\R^n): \left(T\left(\tilde{g}\right)\left(\lambda\right)\right)^\ast\circ T\left(f\right)\left(\lambda\right)v=0 ~~~~~~\, \forall f,g \in L^1\left(\He,\Hi\right)\right\}.$$
    But
    \begin{align*}
        S_\lambda\left(\left( f_1\ast f_2\right)^\lambda\right)
        = S_\lambda\left(f_1^\lambda\ast_\lambda f_2^\lambda\right)= \left(T_\lambda\left(\widetilde{f_2^\lambda}\psi\right)\right)^\ast\circ T_\lambda\left( f_1^\lambda \psi\right)= \left(T_\lambda\left(\Tilde{f_2}\right)\right)^\ast\circ T_\lambda\left( f_1\right),
    \end{align*}
    where in the last equality we have used the fact that $T_\lambda(f)$ depends only on $f^\lambda$ and
    which completes the proof.
\end{proof}

 Note that, if $\rho_\lambda= \pi_\lambda$. Then we have 
 $$\int_{\He}\left(f_1\ast f_2\right)\left(z,t\right)\rho_\lambda\left(z,t\right)\,dz\,dt=  \widehat{\left(\widetilde{f_2}\right)}(\lambda)^\ast\circ \widehat{f_1}(\lambda).$$ So, $\left(T_\lambda\left(\Tilde{f_2}\right)\right)^\ast\circ T_\lambda\left( f_1\right)= \widehat{\left(\widetilde{f_2}\right)}(\lambda)^\ast\circ \widehat{f_1}(\lambda).$ In particular, we have the following

\begin{theorem}
    Let $T: L^1(\mathbb{H}^n,\Hi)\to L^\infty\left(\R^\ast, \mathcal{B}\left(L^2\left(\R^n\right),L^2\left(\R^n,\Hi\right)\right),d\mu\right)$ be a nonzero continuous linear map satisfying: 
    \begin{enumerate}
        \item [(i)] $T\left(\left( f\ast g\right)\fisheye h\right)\left(\lambda\right)=T\left(h\right)\left(\lambda\right)\circ\left( T\left(\tilde{g}\right)\left(\lambda\right)^\ast\circ Tf\left(\lambda\right)\right)$ for all $f,g,h \in L^1(\He,\Hi)$,
        \item [(ii)] $T\left( R_{\left(z,t\right)}f\right)\left(\lambda\right)=T\left(f\right)\left(\lambda\right)\pi_\lambda\left(z,t\right)^\ast,~~~~~~~~ \lambda\in \R^\ast, f\in L^1(\He,\Hi), (z,t)\in \He.$ 
    \end{enumerate}
    Then we have
    $$T\left(g\right)\left(\lambda\right)^\ast \circ Tf\left(\lambda\right)= \left( g\right)\widehat{}\left(\lambda\right)\circ \widehat{f}\left(\lambda\right),~~~~ \lambda\in A, f,g \in L^1(\He,\Hi),$$
    where $A:=\{\lambda\in \R^\ast: Tf(\lambda)\neq 0 \text{ for some } f\in L^1(\He,\Hi)\}$.
\end{theorem}

\begin{proof}
    By previous theorem, for each $\lambda\in A$, there exist a decomposition $L^2\left(\R^n\right)=\mathcal{H}^\lambda\oplus V^\lambda$ and a unitary representation $\rho_\lambda$ on $\He$ such that 
    $$T\left(\tilde{g}\right)\left(\lambda\right)^\ast \circ T\left(f\right)\left(\lambda\right)= \int_{\He}\left(f\ast g\right)\left(w,s\right)\rho_\lambda\left(w,s\right)\,dw\,ds,~~~~~   \text{  on } \mathcal{H}^\lambda,$$
    and $V^\lambda=\left\{v\in L^2(\R^n): T\left(\tilde{g}\right)\left(\lambda\right)^\ast\circ T\left(f\right)\left(\lambda\right)v=0 ~~~~~~,\,\text{for all }~f,g \in L^1\left(\He,\Hi\right)\right\}$. It is enough to show that $\rho_\lambda=\pi_\lambda$. First we show that $V^\lambda=0$.\\
    Using the hypothesis (ii), it follows that $V^\lambda$ is invariant under $\pi_\lambda$. In fact, if $v\in V^\lambda$, then
    \begin{align*}
        \left(T\left(\tilde{g}\right)\left(\lambda\right)^\ast\circ Tf\left(\lambda\right)\right)\left( \pi_\lambda\left(z,t\right)^\ast v\right) = T\left(\tilde{g}\right)\left(\lambda\right)^\ast\circ T\left(R_{\left(z,t\right)}f\right)\left(\lambda\right)= 0.
    \end{align*}
 As $\pi_\lambda$ is irreducible, $V^\lambda= L^2(\R^n)$ or $V^\lambda=\{0\}.$ If $\lambda\in A$, we claim that $V^\lambda=\{0\}.$
In fact, if $V^\lambda= L^2(\R^n),$ then $\left\langle Tf\left(\lambda\right)^\ast Tf\left(\lambda\right)v, v\right\rangle =0$ for all $v\in L^2(\R^n),\,f\in L^1(\He,\Hi)$, which contradicts that $\lambda\in A$.\\
 Now for $\varphi \in L^2(\R^n)$ and $f,g\in L^1(\He,\Hi)$, $\lambda\in A$, consider
    \begin{align*}
        \left(T\left(\tilde{g}\right)\left(\lambda\right)^\ast\circ Tf\left(\lambda\right)\right)\pi_\lambda\left(z,t\right)^\ast \varphi &= T\left(\tilde{g}\right)(\lambda)^\ast \circ \left(Tf\left(\lambda\right)\pi_\lambda\left(z,t\right)^\ast \right)\varphi\\
        &= T\left(\tilde{g}\right)\left(\lambda\right)^\ast\circ T\left(R_{(z,t)}f\right)\left(\lambda\right)\varphi\\
        &= \int_{\He} \left(R_{(z,t)}f\ast g\right)\left(w,s\right)\rho_\lambda\left(w,s\right)\varphi\, dw\, ds\\
        &= \int_{\He} R_{(z,t)}\left(f\ast g\right)\left(w,s\right)\rho_\lambda\left(w,s\right)\varphi\, dw\, ds\\
        &=\int_{\He} \left(f\ast g\right)\left(w,s\right)\rho_\lambda\left(\left(w,s\right).\left(z,t\right)^{-1}\right)\varphi\, dw\, ds\\
        &= \left(\int_{\He} \left(f\ast g\right)\left(w,s\right)\rho\left(w,s\right)\, dw\, ds\right)\circ \rho_\lambda\left(z,t\right)^\ast \varphi\\
        &= \left( T\left(\tilde{g}\right)\left(\lambda\right)^\ast\circ Tf\left(\lambda\right)\right)\circ \rho_\lambda\left(z,t\right)^\ast \varphi.
    \end{align*}
    Thus $$\left( T\left(\tilde{g}\right)\left(\lambda\right)^\ast\circ Tf\left(\lambda\right)\right)\left[\pi_\lambda\left(z,t\right)^\ast\varphi-\rho_\lambda\left(z,t\right)^\ast\varphi\right]=0~~~~~\text{ for all }  f,g\in L^1(\He,\Hi),  $$ which implies that $\pi_\lambda\left(z,t\right)^\ast\varphi-\rho_\lambda\left(z,t\right)^\ast\varphi\in V^\lambda$, for all $\phi \in L^2(\R^n).$ Hence $\pi_\lambda(z,t)=\rho_\lambda(z,t).$ 
\end{proof}

\section{Characterization on $L^2(\He,\Hi)$ }
In this section, we prove the characterization of the Fourier transform on $L^2(\He,\Hi)$. Before stating the result let us observe the action of the Fourier transform on $L^2(\He,\Hi)$.
\par For $f\in L^2(\He)$ and $u\in\Hi$ let $(f\otimes u)(z,t)=f(z,t)u$. Then, for $\phi\in L^2(\R^n)$, we have
$$
(f\otimes u)~\widehat{}~(\lambda)\phi(\xi)=(\widehat{f}(\lambda)\phi)(\xi)\, u=(\widehat{f}(\lambda)\circledast u)\phi(\xi),
$$
where $(T\circledast u)\phi(\xi)= T\phi(\xi)\,u,~\,T\in \mathcal{B}\left(L^2\left(\R^n\right)\right).$ It is immediate that $T\circledast u$ is a bounded linear operator from $L^2(\R^n)$ into $L^2(\R^n,\Hi)$. Moreover, if we let $\mathbf{S}_2$ denote the space of Hilbert-Schmidt operators between the Hilbert spaces $L^2(\R^n)$ and $L^2(\R^n,\Hi)$, then  $\|T\circledast u\|_{\mathbf{S}_2}=\|T\|_{S_2}\|u\|.$ 
\par In other words, $S_2\circledast\Hi$ is isometrically embedded in $\mathbf{S}_2.$ We now show that $S_2\circledast\Hi$ is densely embedded in $\mathbf{S}_2.$ We show that any arbitrary $A\in \mathbf{S}_2$ can be approximated by a finite linear combination of basic tensors. In order to see this let $\{\phi_j:j\in\N\}$ and $\{e_j:j\in\N\}$ be orthonormal bases of $L^2(\R^n)$ and $\Hi$ respectively. Then $\phi_j\otimes e_k$ is an orthonormal basis for $L^2(\R^n,\Hi)$. Choose $a^{(\ell)}_{j,k}\in\C$ such that 
\[
A(e_\ell)=\sum_{j,k\in\N}\,a^{(\ell)}_{j,k}\,\phi_j\otimes e_k.
\]
Then using the fact that $A\in\mathbb{S}_2$, it follows that $\sum_{j,k,\ell\in\N}|a_{j,k,\ell}|^2<\infty$. Now define $T_k$ on $L^2(\R^n)$ by
$$
T_k(\phi_\ell)=\sum_{j\in\N} a^{(\ell)}_{j,k}\phi_j. 
$$
Then it is to show that $\sum_{k\in\N}\|T_k\|^2_{S_2}<\infty$. Further, one can prove that 
$$
\lim_{N\to\infty}\sum_{k=1}^N \left\|T_k\circledast h_k  -A\right\|_{\mathbb{S}_2}=0,
$$
which establishes our claim. With this identification, we have
$
\F_\Hi=\F_\C\boxtimes I
$
on $L^2(\R^n,\Hi)$, where $I$ is the identity operator on $\Hi$ and 
$\mathcal{A}\boxtimes U$ acts on the basic tensors of $\mathbf{S}_2(\equiv \overline{S_2\circledast \Hi}) $ as 
$$
\mathcal{A}\boxtimes U(T\circledast v)= \mathcal{A}(T)\circledast U(v).
$$
Moreover, it follows that the Fourier transform is a map from $L^2(\He,\Hi)$ onto $L^2(\R^\ast,\mathbf{S}_2,d\mu)$.
\par In the remaining part of this section, we aim to characterize the operators on $L^2(\He,\Hi)$ onto $L^2(\R^\ast,\mathbb{S}_2,d\mu)$  which has similar (tensor product) decomposition as the Fourier transform described above. Towards this end we prove the following

\begin{propn} \label{T4.2}
    Let $U:L^2\left( \He,\Hi\right)\to L^2\left(\He,\Hi\right)$ be a bounded linear map of the form $U=S\otimes I$, for some $S\in \mathcal{B}\left(L^2\left(\He\right)\right)$ such that $UR_{(z,t)}f=R_{(z,t)}Uf$ for all $f\in L^2(\He,\Hi)$. Then $\left(Uf\right)\widehat{}\left(\lambda\right)=\left(m\left(\lambda\right)\boxtimes I\right)\widehat{f}\left(\lambda\right)$, for some $m\in L^\infty\left(\R^\ast, \mathcal{B}\left(L^2\left(\R^n\right)\right),d\mu\right).$
\end{propn}

\begin{proof}
    By making use of the fact that $U$ commutes with right translations, it follows that $S$ also commutes with the right translations. In fact, if $f=F\otimes u$ is a simple tensor then
    $$
    R_{(z,t)}SF\otimes u= UR_{(z,t)}f= UR_{(z,t)}f = SR_{(z,t)}F\otimes u,
    $$
    for all $u\in\Hi$. Therefore, $SR_{(z,t)}= R_{(z,t)}S.$
Hence, by \cite{hormander1960estimates, Lavanya_Heisen}, there exists $m\in L^\infty\left(\R^\ast, \mathcal{B}\left(L^2\left(\R^n\right)\right),d\mu\right)$ such that $\left( SF\right)\widehat{}\left(\lambda\right)= m\left(\lambda\right)\widehat{f}\left(\lambda\right).$ Thus
    \begin{align*}
        \left( Uf\right)\widehat{}\left(\lambda\right)
        &= \left( SF\right)\widehat{}\left(\lambda\right)\otimes u= m\left(\lambda\right)\widehat{F}\left(\lambda\right)\otimes u= \left( m\left(\lambda\right)\boxtimes I\right)\widehat{f}\left(\lambda\right).
    \end{align*}
\end{proof}
Now we are in a position to prove the main theorem of this section. We have the following characterization of the Fourier transform on $L^2(\He,\Hi)$.
\begin{theorem} \label{T2.2}
    Let $T: L^2(\He,\Hi)\to L^2\left(\R^\ast, \mathbf{S}_2,d\mu\right)$ be a nonzero surjective linear map satisfying:
    \begin{enumerate}
        \item [(i)] $T\left(\left( f\ast g\right)\fisheye h\right)\left(\lambda\right)=Th\left(\lambda \right)\circ \left(T\tilde{g}\left(\lambda\right)^\ast \circ Tf\left(\lambda\right)\right)$ for all $f,g,h\in L^2(\He,\Hi)$,
        \item [(ii)] $T\left(R_{(z,t)}f\right)\left(\lambda\right)=Tf\left(\lambda\right)\pi_\lambda\left(z,t\right)^\ast$, $\lambda\in \R^\ast$, $f\in L^2(\He,\Hi)$, $(z,t)\in \He$.
        \item [(iii)] $T=S\boxtimes I$, where $S$ is a bounded linear operator from $L^2(\He)$ into $L^2\left(\R^\ast, S_2(L^2(\R^n)), d\mu\right)$.
    \end{enumerate}
    Then $$Tg\left(\lambda\right)^\ast\circ Tf\left(\lambda\right)= \widehat{g}\left(\lambda\right)\circ \widehat{f}\left(\lambda\right)~~~~\text{ for all } f,g\in L^2(\He,\Hi), \lambda\in \R^\ast. $$
\end{theorem}

\begin{proof}
    Define $U: L^2(\He,\Hi)\to L^2(\He,\Hi)$ by $Uf=g$ if $Tf=\widehat{g}$. Thus $\left(Uf\right)\widehat{}\left(\lambda\right)=Tf\left(\lambda\right).$ Clearly, $U$ is surjective.
    \par Let $f_1,f_2,f_3\in L^2(\He,\Hi)$ be such that $\left(f_1\ast f_2\right)\fisheye f_3\in L^2(\He,\Hi)$. Then
    \begin{align*}
        \left(U\left(\left(f_1\ast f_2\right)\fisheye f_3\right)\right)~\widehat{}~\left(\lambda\right)
        &= T\left(\left(f_1\ast f_2\right)\fisheye f_3\right)\left(\lambda\right)\\
        &= Tf_3\left(\lambda\right)\circ \left( T\tilde{f_2}\left(\lambda\right)^\ast \circ Tf_1\left(\lambda\right)\right)\\
        &= \left(Uf_3\right)~\widehat{}~\left(\lambda\right)\circ \left(\left( \left( U\tilde{f_2}\right)~\widehat{}~\left(\lambda\right)\right)^\ast\circ \left(Uf_1\right)~\widehat{}~\left(\lambda\right)\right)\\
        &= \left(\left(Uf_1\ast \widetilde{U\tilde{f_2}}\right)\fisheye Uf_3\right)~\widehat{}~\left(\lambda\right),
    \end{align*}
    which implies that $U\left( \left( f_1\ast f_2\right)\fisheye f_3\right)= \left(Uf_1\ast \widetilde{U\tilde{f_2}}\right)\fisheye Uf_3.$
    Further, $U$ commutes with the right translations. In fact, 
    \begin{align*}
        \left( UR_{(z,t)}f\right)~\widehat{}~\left(\lambda \right)
        = T\left( R_{(z,t)}f\right)\left(\lambda\right)= Tf\left(\lambda\right)\pi_\lambda\left(z,t\right)^\ast= \left(Uf\right)~\widehat{}~\left(\lambda\right)\pi_\lambda\left(z,t\right)^\ast=\left(R_{(z,t)}Uf\right)~\widehat{}~\left(\lambda\right).
    \end{align*}
    Furthermore, we claim that $U=S_1\otimes I$ for some $S_1\in \mathcal{B}\left(L^2(\He) \right).$ In order to see this,  define $S_1: L^2(\He)\to L^2\left( \He\right)$ by $S_1(F)= G$ if $S(F)=\widehat{G}.$ Then using the fact that $T=S\boxtimes I$, it follows that $U$ agrees with $S_1\otimes I$ on set of all basic tensors. Hence $U= S\otimes I$ on $L^2(\He,\Hi).$ Then, appealing to Proposition \ref{T4.2}, there exists $m\in L^\infty\left(\R^\ast, \mathcal{B}\left(L^2\left(\R^n\right)\right),d\mu\right)$ such that 
    $$
    \left(Uf\right)~\widehat{}~\left(\lambda\right)=\left(m\left(\lambda\right)\boxtimes I\right)\widehat{f}\left(\lambda\right).
    $$
    Now for $f,g,h\in L^2(\He,\Hi)$ with $\left(f\ast g\right)\fisheye h\in L^2(\He,\Hi)$ we have $$U\left(\left(f\ast g\right)\fisheye h\right)~\widehat{}~\left(\lambda\right)= \left(m\left(\lambda\right)\boxtimes I\right)\left(\left( f\ast g\right)\fisheye h\right)~\widehat{}~\left(\lambda\right).$$
    But
    \begin{align*}
        U\left(\left(f\ast g\right)\fisheye h\right)~\widehat{}~\left(\lambda\right) &= \left(\left( Uf\ast \widetilde{U\tilde{g}}\right)\fisheye Uh\right)~\widehat{}~\left(\lambda\right)\\
        &= \left(Uh\right)~\widehat{}~\left(\lambda\right)\circ \left( \left( \widetilde{\widetilde{U\tilde{g}}}\right)~\widehat{}~\left(\lambda\right)^\ast\circ \left(Uf\right)~\widehat{}~\left(\lambda\right)\right)\\
        &= \left(Uh\right)~\widehat{}~\left(\lambda\right)\circ \left( \left( U\tilde{g}\right)~\widehat{}~\left(\lambda\right)^\ast\circ \left(Uf\right)~\widehat{}~\left(\lambda\right)\right),
    \end{align*}
    and
    \begin{align*}
        \left(m\left(\lambda\right)\boxtimes I\right)\left(\left( f\ast g\right)\fisheye h\right)~\widehat{}~\left(\lambda\right)
        &=\left( m\left(\lambda \right)\boxtimes I\right)\left( \widehat{h}\left(\lambda\right)\circ\left(\widehat{\tilde{g}}\left(\lambda\right)^\ast \circ \widehat{f}\left(\lambda\right)\right)\right)\\
        &=\left(\left( m\left(\lambda \right)\boxtimes I\right) \widehat{h}\left(\lambda\right)\right)\circ\left(\widehat{\tilde{g}}\left(\lambda\right)^\ast \circ \widehat{f}\left(\lambda\right)\right)\\
        &= \left(Uh\right)~\widehat{}~\left(\lambda\right)\circ \left(\widehat{\tilde{g}}\left(\lambda\right)^\ast \circ \widehat{f}\left(\lambda\right)\right).
    \end{align*}
    Therefore,
    $$\left(Uh\right)~\widehat{}~(\lambda)\circ \left( \widehat{\tilde{g}}\left(\lambda\right)^\ast\circ \widehat{f}\left(\lambda\right)- \left(U\tilde{g}\right)~\widehat{}~\left(\lambda\right)^\ast\circ \left(Uf\right)~\widehat{}~\left(\lambda\right)\right)=0,~~\text{ for all } h\in L^2(\He, \Hi).$$
    Therefore by surjectivity of $U$, we have
    $$\widehat{h}(\lambda)\circ \left( \widehat{\tilde{g}}\left(\lambda\right)^\ast\circ \widehat{f}\left(\lambda\right)- \left(U\tilde{g}\right)~\widehat{}~\left(\lambda\right)^\ast\circ \left(Uf\right)~\widehat{}~\left(\lambda\right)\right)=0,$$ 
    which implies that
    $$Range\left(\widehat{\tilde{g}}\left(\lambda\right)^\ast\circ \widehat{f}\left(\lambda\right)- \left(U\tilde{g}\right)~\widehat{}~\left(\lambda\right)^\ast\circ \left(Uf\right)~\widehat{}~\left(\lambda\right)\right)\subseteq Ker\left( \widehat{h}(\lambda)\right), \text{ for all } h\in L^2(\He,\Hi).$$ 
    This forces $\widehat{\tilde{g}}\left(\lambda\right)^\ast\circ \widehat{f}\left(\lambda\right)- \left(U\tilde{g}\right)~\widehat{}~\left(\lambda\right)^\ast\circ \left(Uf\right)~\widehat{}~\left(\lambda\right)=0,$ 
    which is same as $T\tilde{g}\left(\lambda\right)^\ast\circ Tf\left(\lambda\right)=\widehat{\tilde{g}}\left(\lambda\right)^\ast\circ \widehat{f}\left(\lambda\right).$ Therefore,
     
    $$Tg\left(\lambda\right)^\ast\circ Tf\left(\lambda\right)=\widehat{g}\left(\lambda\right)^\ast\circ \widehat{f}\left(\lambda\right) \text{ for all } f,g \in L^2(\He,\Hi).$$
\end{proof}

\section{Characterization on Schwartz space}
It is well known that the Schr\"{o}dinger representations $\pi_\lambda$  are all infinite-dimensional irreducible unitary representations of $\He.$ In addition, $\He$ admits a family of certain one-dimensional representations. Indeed, for each $\zeta\in \C^n$, the character $\chi_\zeta(z,t)= e^{i\Re z\cdot\overline{\zeta}}$ defines those one dimensional representations of $\He.$ Then, using the representations, 
$$\pi_{\lambda, \zeta}(z,t)= \chi_\zeta(z,t)\pi_\lambda(z,t)$$
parametrized by $(\lambda, \zeta)\in \R^\ast\times \C^n,$ a modified version of the group Fourier transform, known as the ``Fourier–Weyl transform", is defined in \cite{Lavanya_Revisiting}. Furthermore, they established a convolution theorem for convolution in the central variable and subsequently used this result, together with the ordinary convolution theorem, to obtain a characterization of the Fourier–Weyl transform. As a consequence, they derived a characterization of the group Fourier transform. In this section, our aim is to establish a vector-valued analogue of this result.

Following \cite{Lavanya_Revisiting}, we define the Fourier-Weyl transform for $f\in L^1(\He,\Hi)$ as
$$\F f(\lambda, \zeta):= \int_{\He}f(z,t)\pi_{\lambda,\zeta}(z,t)\, dz\, dt,$$
where the above integral is a Bochner integral, see \cite{thangavelu1993paley}.

Recall the definition of the Group Fourier transform of a function $f$ can be given by
$$\widehat{f}(\lambda)= \int_{\He} f(z,t)\pi_\lambda(z,t)\, dz\, dt.$$
By making use of the following easily verifiable equality, 
$$\pi_{\lambda, \zeta}\left(z,t\right)= \pi_\lambda\left(-i\lambda^{-1}\zeta,0\right)\pi_\lambda\left(z,t\right)\pi_\lambda\left(i\lambda^{-1}\zeta,0\right),$$
we have
$$\F f\left(\lambda,\zeta\right)=\pi_\lambda\left(-i\lambda^{-1}\zeta,0\right)\widehat{f}(\lambda) \pi_\lambda\left(i\lambda^{-1}\zeta,0\right).$$
By analogy with this identity, we introduce the following notation for any $S\in L^2(\R^\ast, S_2,d\mu)$ in the sequel.
$$
S\left(\lambda,\zeta\right)=\pi_\lambda\left(-i\lambda^{-1}\zeta,0\right)S(\lambda)\pi_\lambda\left(i\lambda^{-1}\zeta,0\right).
$$
Motivated by the work of \cite{Lavanya_Revisiting}, given $f,g \in L^1(\He,\Hi)$, we define $f\ast_3 g$, the $*$-convolution in the central variable by
$$f\ast_3 g(z,t)=\int_{-\infty}^\infty \left\langle f\left( z,t-s\right), g\left(z,s\right)\right\rangle \, ds.$$
Similarly for $\phi\in L^1(\He)$ and $h\in L^1(\He,\Hi)$, we define $\phi\fisheye_3 h$, the $\fisheye$-convolution in the central variable by
$$\phi \fisheye_3 h(z,t)= \int_{-\infty}^\infty \phi(z,t-s)h(z,s)\, ds.$$

Now we observe how Fourier-Weyl transform behaves with these convolutions.
\begin{propn}
    Let $f,g\in L^1(\He,\Hi)$, then we have the following:
    $$\F (f\ast_3 g)(\lambda,z)=\int_{\C^n}\left\langle \F f\left(\lambda,\omega\right)\left(\cdot\right), \f_{\C^n\times \R}^{-1}\tilde{g}\left(z-w,\lambda\right)\right\rangle\, dw,$$

where the above is interpreted in the weak sense: for $\phi\in L^2(\R^n)$,
$$\F (f\ast_3 g)(\lambda,z)\phi(x)=\int_{\C^n}\left\langle \F f\left(\lambda,\omega\right)\phi\left(x\right), \f_{\C^n\times \R}^{-1}\tilde{g}\left(z-w,\lambda\right)\right\rangle\, dw.$$
Similarly, if $\phi\in L^1(\He)$ and $h\in L^1(\He,\Hi)$, we have 
$$\F \left( \phi\fisheye_3 h\right)\left(\lambda,z\right)= \int_{\C^n} \F \phi(\lambda, w)(\cdot)\f_{\C^n\times \R}^{-1} h(z-w,\lambda)\, dw.$$
\end{propn}

\begin{proof}
Let $h\in L^1(\He,\C)$ and $\phi,\,\psi\in L^2(\R^n)$. Then from the definition of the Fourier-Weyl transform, it follows that
\begin{equation}\label{Eq 5.1}
    \left( \F h(\lambda, z)\phi, \psi\right)=\int_{\C^n} h^\lambda(w)e^{i\Re(w\overline{z})}\left(\pi_\lambda(w)\phi, \psi\right)\, dw,
\end{equation} 
which is same as
\begin{equation}\label{Eq 5.2}
    \f_{\C^n}^{-1}\left( h^\lambda(\cdot)\left( \pi_\lambda(\cdot)\phi, \psi\right)\right)(z)=\left(\F h(\lambda, z)\phi, \psi\right).
\end{equation}
Here the notation $(f,g)$ for $f\in L^2(\R^n,\Hi)$, $g\in L^2(\R^n)$ stands as follows 
$$(f,g):= \int_{\R^n}\overline{g(x)}f(x)\, dx.$$
    Let $f,g\in L^1(\He,\Hi)$. Then $f\ast_3 g\in L^1(\He,\C)$. Thus, using \ref{Eq 5.1}, \ref{Eq 5.2} and the convolution theorem for the Fourier transform of vector valued functions, we have
    \begin{align*}
        \left(\F \left(f\ast_3 g\right)\left(\lambda, z\right)\phi,\psi\right)
        &= \int_{\C^n} \left( f\ast_3 g\right)^\lambda\left(w\right)e^{i\Re\left(w\overline{z}\right)}\left(\pi_\lambda\left(w\right)\phi, \psi\right)\, dw\\
        &= \int_{\C^n} \left\langle \left(\pi_\lambda\left(w\right)\phi, \psi\right)f^\lambda\left(w\right), g^{-\lambda}\left(w\right)\right\rangle e^{i\Re\left(w\overline{z}\right)}\, dw\\
        &= \int_{\C^n} \left\langle \f_{\C^n}^{-1}\left( \left( \pi_\lambda\left(\cdot\right)\phi,\psi\right)f^\lambda\left(\cdot\right)\right)\left(w\right), \f_{\C^n}^{-1}\tilde{g}^\lambda\left(z-w\right)\right\rangle\, dw\\
        &= \int_{\C^n} \left\langle \left( \F f\left(\lambda, w\right)\phi, \psi\right), \f_{\C^n\times \R}^{-1}\tilde{g}\left(z-w,\lambda\right)\right\rangle \, dw\\
        &= \int_{\C^n}\left\langle\int_{\R} \overline{\psi\left(x\right)}\F f\left(\lambda,w\right)\phi\left(x\right)\, dx, \f_{\C^n\times\R}^{-1}\tilde{g}\left(z-w, \lambda\right)\right\rangle\, dw\\
         &=\left(\left(\int_{\C^n} \left\langle\F f\left(\lambda,w\right)\left(\cdot\right), \f_{\C^n\times\R}^{-1}\tilde{g}\left(z-w, \lambda\right)\right\rangle\, dw\right)\left(\phi\right),\psi\right).  
    \end{align*}
    Since the above equality holds for any $\phi,\psi$, we get
    $$\F \left(f\ast_3 g\right)\left(\lambda,z\right)= \int_{\C^n} \left\langle\F f\left(\lambda,w\right)\left(\cdot\right), \f_{\C^n\times\R}^{-1}\tilde{g}\left(z-w, \lambda\right)\right\rangle\, dw. $$
\end{proof}
As a consequence of the above proposition, for $f,g,h\in L^1(\He,\Hi)$, we get
\begin{align} \label{E5.3}
    \F\left(\left( f\ast_3 g\right)\fisheye_3 h\right)\left( \lambda, z\right)
    &= \int_{\C^n} \F\left( f\ast_3 g\right)\left( \lambda, w\right)\left(\cdot\right) \f_{\C^n\times \R}^{-1} h\left( z-w, \lambda\right)\, dw \nonumber\\
    &= \int_{\C^n}\int_{\C^n} \left\langle \F f\left(\lambda, \zeta\right)\left(\cdot\right), \f_{\C^n\times \R}^{-1}\tilde{g}\left( w-\zeta,\lambda\right)\right\rangle \f_{\C^n\times \R}^{-1}h\left( z-w,\lambda\right)\, d\zeta\, dw.
\end{align}
We express the right-hand side of the above equality in an alternative form as follows.
\begin{propn}
    For $f,g,h\in L^1(\He, \Hi)$, we have
    \begin{equation}\label{eq: cent conv alt}
    \F\left(\left(f\ast_3 g\right)\fisheye_3 h\right)(\lambda, \zeta)= \int_{\C^n} \f_{\C^n\times\R}^{-1}f \ast \f_{\C^n\times \R}^{-1}\widetilde{g}(\zeta-z, \lambda)\widehat{h}(\lambda,z)\, dz.
    \end{equation}
\end{propn}

\begin{proof}
    It is easy to see that 
    \begin{equation}\label{Eq 5.4}
    \left(\left(f\ast_3 g\right)\fisheye_3 h\right)^\lambda(z)= \left\langle f^\lambda(z), g^{-\lambda}(z)\right\rangle h^\lambda(z).
    \end{equation}
    Now, for $\phi\in L^2(\R^n),\,\psi\in L^2(\R^n,\Hi)$, using \ref{Eq 5.1},\ref{Eq 5.2} and \ref{Eq 5.4} we have,
    \begin{align*}
        \left( \F\left(\left(f\ast_3 g\right)\fisheye_3 h\right)(\lambda, \zeta)\phi, \psi\right)&= \int_{\C^n} \left\langle f^\lambda(w), g^{-\lambda}(w)\right\rangle h^\lambda(w)e^{i\Re(w\overline{\zeta})}\left(\pi_\lambda(w)\phi, \psi\right)\, dw\\
        &= \int_{\C^n} \left(\left\langle f^\lambda(w), g^{-\lambda}(w)\right\rangle\right)\left(h^\lambda(w)\left(\pi_\lambda(w)\phi, \psi\right)\right)e^{i\Re(w\overline{\zeta})}\, dw\\
        &= \f_{\C^n}^{-1}\left( \left\langle f^\lambda(\cdot), g^{-\lambda}(\cdot)\right\rangle h^\lambda(\cdot)\left(\pi_\lambda(\cdot)\phi, \psi\right)\right)(\zeta)\\
        &= \f_{\C^n}^{-1}\left(\left\langle f^\lambda(\cdot), g^{-\lambda}(\cdot)\right\rangle\right)\fisheye \f_{\C^n}^{-1}\left(h^\lambda(\cdot)\left(\pi_\lambda(\cdot)\phi, \psi\right)\right)(\zeta)\\
        &= \int_{\C^n} \f_{\C^n}^{-1}\left(\left\langle f^\lambda(\cdot), g^{-\lambda}(\cdot)\right\rangle\right)(\zeta-z)\f_{\C^n}^{-1}\left( h^\lambda(\cdot)\left( \pi_\lambda(\cdot)\phi,\psi\right)\right)(z)\, dz\\
        &=  \int_{\C^n} \f_{\C^n}^{-1}\left(\left\langle f^\lambda(\cdot), g^{-\lambda}(\cdot)\right\rangle\right)(\zeta-z)\left( \F h(\lambda,z)\phi,\psi\right)\, dz\\
        &= \int_{\C^n} \f_{\C^n}^{-1}f^\lambda\ast \f_{\C^n}^{-1}\widetilde{g^{-\lambda}}(\zeta-z)\left( \F h(\lambda,z)\phi, \psi\right)\, dz\\
        &= \int_{\C^n}\f_{\C^n\times \R}^{-1}f \ast \f_{\C^n\times \R}^{-1}\widetilde{g}(\zeta-z,\lambda)\left(\F h(\lambda, z)\phi, \psi\right)\, dz
    \end{align*}
    As above equality is true for any $\phi ,\psi$,
    $$\F\left(\left(f\ast_3 g\right)\fisheye_3 h\right)(\lambda, \zeta)= \int_{\C^n} \f_{\C^n\times\R}^{-1}f \ast \f_{\C^n\times \R}^{-1}\widetilde{g}(\zeta-z, \lambda)\F h(\lambda,z)\, dz.$$
\end{proof}

\begin{remark} In the following theorem, we characterize the operators $T$ that share properties similar to \eqref{E5.3} and \eqref{eq: cent conv alt}, which are satisfied by the Fourier–Weyl transform. We denote $\pi_\lambda(-i\lambda^{-1}\zeta,0)Tf(\lambda)\pi_\lambda(i\lambda^{-1}\zeta,0)$ by $Tf(\lambda,\zeta)$. Then $Tf(\lambda,0)= Tf(\lambda).$ Thus, we consider operators $T$ satisfying the following conditions:  
\begin{align}
\label{eq: prop 1}T\left(\left( f\ast_3 g\right)\fisheye_3 h\right)\left(\lambda\right)&=
\int_{\C^n}\int_{\C^n} \left\langle Tf\left(\lambda, \zeta\right)\left(\cdot\right), \f_{\C^n\times \R}^{-1}\tilde{g}\left( w-\zeta,\lambda\right)\right\rangle \f_{\C^n\times \R}^{-1}h\left(-w,\lambda\right)\, d\zeta\, dw.\\
\label{eq: prop 2}T\left(\left( f\ast_3 g\right)\fisheye_3 h\right)\left(\lambda\right)&= \int_{\C^n} \f_{\C^n\times\R}^{-1}f \ast \f_{\C^n\times \R}^{-1}\widetilde{g}(-z, \lambda)Th(\lambda,z)\, dz.    
\end{align}
\end{remark}

Let $S(\He,\Hi)$ denote the Schwartz space on $\He$, taking values in a separable Hilbert space $\Hi,$ and $\widehat{S}(\He)$ denote
the image of $S(\He, \Hi)$ under the group Fourier transform $f\to \widehat{f}(\lambda).$ We can view $\widehat{S}(\He,\Hi)$ as a subspace, embedded in $L^2(\R^\ast, \mathbf{S}_2, d\mu).$ Here $\mathbf{S}_2$ is the space of Hilbert-Schmidt operators between the Hilbert spaces $L^2(\R^n)$ and $L^2(\R^n,\Hi)$. 

\begin{theorem}\label{T5.2}
Let $T:S(\He,\Hi)\to \widehat{S}(\He,\Hi)$ be a bijection satisfying \eqref{eq: prop 1}, \eqref{eq: prop 2} and
\begin{equation}
T\left( \left( f\ast g\right)\fisheye h\right)\left(\lambda\right)= Th(\lambda)\circ \left( T\Tilde{g}\left(\lambda\right)^\ast \circ Tf\left(\lambda\right)\right),~~~~~~~
\end{equation}
for all $f,g,h \in S(\He, \Hi )$.
Then there exist functions $\alpha,\,\beta: \R^\ast\to \C^n$ such that $$Tf(\lambda)=\alpha(\lambda) \F f\left(\lambda, \beta\left(\lambda\right)\right),~~\text{  for all }f\in S(\He,\Hi).$$
\end{theorem}
\begin{proof}
     Define $U:S(\He,\Hi)\to S(\He,\Hi)$ by $Uf=g$, if $Tf(\lambda)=\widehat{g}(\lambda)$. Clearly, $U$ is a bijection and $Tf(\lambda)= \widehat{Uf}(\lambda)$, $Tf(\lambda,\zeta)= \F (Uf)(\lambda,\zeta)$ for all $f\in S(\He,\Hi).$

    Further, for $f,g,h\in S(\He,\Hi)$, using the hypothesis (i) we have
    \begin{align*}
        \left(U\left(\left( f\ast_3 g\right)\fisheye_3 h\right)\right)~\widehat{}~(\lambda) &= T\left(\left( f\ast_3 g\right)\fisheye_3 h\right)\left(\lambda\right)\\
        &= \int_{\C^n}\int_{\C^n} \left\langle Tf\left(\lambda, \zeta\right)\left(\cdot\right), \f_{\C^n\times\R}^{-1}\tilde{g}\left(w-\zeta,\lambda\right)\right\rangle\f_{\C^n\times \R}^{-1}h\left(-w,\lambda\right)\,dw\, d\zeta\\
        &= \int_{\C^n}\int_{\C^n} \left\langle \F\left(Uf\right)\left(\lambda, \zeta\right)\left(\cdot\right), \f_{\C^n\times\R}^{-1}\tilde{g}\left(w-\zeta,\lambda\right)\right\rangle\f_{\C^n\times \R}^{-1}h\left(0-w,\lambda\right)\,dw\, d\zeta\\
        &= \F\left(\left( Uf\ast_3 g\right)\fisheye_3 h\right)\left(\lambda,0\right)\\
        &= \left(\left( Uf\ast_3 g\right)\fisheye_3 h\right)~\widehat{}~(\lambda).
    \end{align*}
    Therefore we have, $$U\left(\left( f\ast_3 g\right)\fisheye_3 h\right)= \left( Uf\ast_3 g\right)\fisheye_3 h.$$

 In a similar manner using hypothesis (ii) we can show that 
 $$U\left(\left( f\ast_3 g\right)\fisheye_3 h\right)= \left(f\ast_3 g\right) \fisheye_3 Uh.$$
 Hence  $\left( Uf\ast_3 g\right)\fisheye_3 h= \left( f\ast_3 g\right)\fisheye_3 Uh$. By taking the inverse Fourier transform in the central variable on both sides, it follow that 
 $\left\langle \left(Uf\right)^\lambda\left(z\right), g^{-\lambda}\left(z\right)\right\rangle h^\lambda\left(z\right)= \left\langle f^\lambda\left(z\right), g^{-\lambda}\left(z\right)\right\rangle \left(Uh\right)^\lambda\left(z\right) $ for every $\lambda\in \R^\ast.$\\
 \par Now, for a given $\lambda\in \R^\ast,$ choose $f,g\in S(\He,\Hi)$ be such that $\left\langle f^\lambda(z), g^{-\lambda}(z)\right\rangle$ vanishes nowhere on $\C^n.$ Then for any $h\in S(\He,H)$
  $$\left(Uh\right)^\lambda\left(z\right)= \frac{\left\langle \left(Uf\right)^\lambda\left(z\right), g^{-\lambda}\left(z\right)\right\rangle}{\left\langle f^\lambda\left(z\right), g^{-\lambda}\left(z\right)\right\rangle}h^\lambda\left(z\right)= \varphi_\lambda(z)h^\lambda(z),$$
  where
  $$
  \varphi_\lambda(z)=\frac{\left\langle \left(Uf\right)^\lambda\left(z\right), g^{-\lambda}\left(z\right)\right\rangle}{\left\langle f^\lambda\left(z\right), g^{-\lambda}\left(z\right)\right\rangle}
  $$
  which is clearly continuous.
Furthermore, using the hypothesis (iii) we obtain
    \begin{align*}
        \left(U\left(\left( f\ast g\right)\fisheye h\right)\right)~\widehat{}~(\lambda) &= T\left(\left( f\ast g\right)\fisheye h\right)\left(\lambda\right)\\
        &= Th\left(\lambda\right)\circ \left( T\tilde{g}\left(\lambda\right)^\ast \circ Tf\left(\lambda\right)\right)\\
        &= \left(Uh\right)~\widehat{}~(\lambda)\circ \left(\left(U\tilde{g}\right)~\widehat{}~(\lambda)^\ast\circ\left(Uf\right)~\widehat{}~(\lambda)\right)\\
        &= \left(\left(Uf\ast \widetilde{U\widetilde{g}}\right)\fisheye Uh\right)~\widehat{}~(\lambda),
    \end{align*}
    which implies that $$U\left(\left(f\ast g\right)\fisheye h\right)= \left( Uf\ast \widetilde{U\widetilde{g}}\right)\fisheye Uh,\,\,\,\, f,g,h\in S(\He,\Hi).$$
{\color{black}
Hence, for $f,g,h \in S(\He,\Hi)$
$$
\varphi_\lambda(z)\left(\left(f\ast g\right)\fisheye h\right)^\lambda\left(z\right)=\left(U\left(\left(f\ast g\right)\fisheye h\right)\right)^\lambda \left(z\right)= \left(\left( Uf\ast \widetilde{U\widetilde{g}}\right)\fisheye Uh\right)^\lambda\left(z\right).
$$
We simplify the right hand side as
\begin{align*}
&\left(\left(Uf\right)^\lambda\ast_\lambda\left(\widetilde{U\widetilde{g}}\right)^{-\lambda}\right)\fisheye_\lambda\left(Uh\right)^\lambda(z)=\\
& \int_{\C^n}\int_{\C^n} \varphi_\lambda(x)\overline{\varphi_\lambda(x-w)}\varphi_\lambda(z-w)\left\langle f^\lambda(x), g^{-\lambda}(w-x)\right\rangle h^\lambda(z-w) e^{i\frac{\lambda}{2}\Im(z\overline{w})}e^{i\frac{\lambda}{2}\Im(w\overline{x})}\, dw\, dx
\end{align*}

On the other hand,
\begin{align*}
    \varphi_\lambda(z)\left(\left(f\ast g\right)\fisheye h\right)^\lambda\left(z\right)
    &= \varphi_\lambda(z)\left(\left(f^\lambda\ast_\lambda g^{-\lambda}\right)\fisheye_\lambda h^\lambda\right)^\lambda\left(z\right)\\
    &= \int_{\C^n}\int_{\C^n} \varphi_\lambda(z)\left\langle f^\lambda(x), g^{-\lambda}(w-x)\right\rangle h^\lambda(z-w) e^{i\frac{\lambda}{2}\Im(z\overline{w})}e^{i\frac{\lambda}{2}\Im(w\overline{x})}\, dw\, dx.
\end{align*}
Therefore 
$$\int_{\C^n}\int_{\C^n} \left(\varphi_\lambda(z)-\varphi_\lambda(x)\overline{\varphi_\lambda(x-w)}\varphi_\lambda(z-w)\right)\left\langle f^\lambda(x), g^{-\lambda}(w-x)\right\rangle h^\lambda(z-w) e^{i\frac{\lambda}{2}\Im(z\overline{w})}e^{i\frac{\lambda}{2}\Im(w\overline{x})}\, dw\, dx=0,$$
for all $f,g,h \in S(\He, \Hi).$
}

By choosing appropriate $f,g,h\in S(\He,\Hi)$, it is easy to see that 
\begin{equation}\label{eq: phi prop}
    \varphi_\lambda(z)=\varphi_\lambda(x)\overline{\varphi_\lambda(w-x)}\varphi_\lambda(z-w) \text{ for all }x,w,z\in \C^n.
\end{equation} 
Moreover, $|\varphi_\lambda(0)|=1$ and $\varphi_\lambda(z)\neq 0$ for all $z\in \C^n.$ In fact, if $\varphi_\lambda(z_0)=0$, then 
$$\varphi_\lambda(z)= \varphi_\lambda(z_0)\overline{\varphi_\lambda(w-z_0)}\varphi_\lambda(z-w)=0,~~\text{ for all }\, z\in \C^n.$$
Again, by putting $x=z,~ w=0$ in \eqref{eq: phi prop} we get $\overline{\varphi_\lambda(-z)}=\frac{1}{\varphi_\lambda(z)}$. Thus, we can rewrite \eqref{eq: phi prop}, as 
$$\varphi_\lambda(z)\varphi_\lambda(x-w)=\varphi_\lambda(x)\varphi_\lambda(z-w).$$
Now taking $z=z+w$ and then $x=w$, we get
$$\frac{\varphi_\lambda(z+w)}{\varphi_\lambda(0)}=\frac{\varphi_\lambda(z)}{\varphi_\lambda(0)}\cdot\frac{\varphi_\lambda(w)}{\varphi_\lambda(0)}.$$
Let $\eta_\lambda(z)=\frac{\varphi_\lambda(z)}{\varphi_\lambda(0)}.$ Then $\eta_\lambda(z+w)=\eta_\lambda(z)\eta_\lambda(w)$ for all $z,w\in \C^n.$ 
\par Since $U$ is a bijection on $S(\He,\Hi)$ and $\left(Uf\right)^\lambda\left(z\right)=\varphi_\lambda\left(z\right)f^\lambda\left(z\right), f\in S(\He,\Hi)$, it is immediate that $\varphi_\lambda$ (hence $\eta_\lambda$) is bounded. Thus, $\varphi_\lambda(z)=\varphi_\lambda(0)e^{-i\Re\left(z.\overline{\beta\left(\lambda\right)}\right)}$ for some $\beta(\lambda)\in \C^n$.

\par
  Now for $f\in S(\He,\Hi)$, we have
  \begin{align*}
      Tf(\lambda)&= (Uf)~\widehat{}~(\lambda)= \int_{\C^n}\left(Uf\right)^\lambda\left(z\right)\pi_\lambda\left(z\right)\, dz
      = \varphi_\lambda(0)\int_{\C^n} e^{-i\Re\left(z.\overline{\beta\left(\lambda\right)}\right)}f^\lambda\left(z\right)\pi_\lambda\left(z\right)\, dz.
  \end{align*}
  Hence, $Tf(\lambda)= \varphi_\lambda(0)\F f\left(\lambda, \beta\left(\lambda\right)\right),~~~~~ f\in S(\He,\Hi),$ which completes the proof.
\end{proof}
\begin{corollary}
    In addition to the hypothesis of Theorem \ref{T5.2}, if $T$ satisfies 
    \begin{equation}
        T\left(R_{\left(z,0\right)}f\right)\left(\lambda\right)=Tf\left(\lambda\right)\pi_\lambda\left(z\right)^\ast,
    \end{equation}
    then $Tf(\lambda)=\varphi_\lambda(0)\widehat{f}(\lambda).$
\end{corollary}
\begin{proof}
    From the above hypothesis we have $T\left(R_{\left(z,0\right)}f\right)\left(\lambda\right)=Tf\left(\lambda\right)\pi_\lambda\left(z\right)^\ast$ and by previous theorem we get $Tf(\lambda)= \varphi_\lambda(0)\F f(\lambda,\beta(\lambda))$. Thus 
    $$
    Tf(\lambda)\pi_\lambda(z)^\ast= \varphi_\lambda(0)\F\left(R_{(z,0)}f\right)\left(\lambda,\beta\left(\lambda\right)\right),
    $$
    which is same as
    $$
    \varphi_\lambda(0)\F f\left(\lambda,\beta(\lambda)\right)\pi_\lambda(z)^\ast =\varphi_\lambda(0)\F\left(R_{(z,0)}f\right)\left(\lambda,\beta\left(\lambda\right)\right).
    $$
    Upon simplification we arrive at
    $$
    \widehat{f}(\lambda)\pi_\lambda\left(i\lambda^{-1}\beta\left(\lambda\right)\right)\pi_\lambda(z)^\ast = \widehat{f}(\lambda)\pi_\lambda\left(z\right)^\ast\pi_\lambda\left(i\lambda^{-1}\beta\left(\lambda\right)\right), ~~\text{ for all } z\in\C^n,\,f\in S(\He,\Hi),
    $$
    which implies that
    $$\pi_\lambda\left(i\lambda^{-1}\beta\left(\lambda\right)\right)\pi_\lambda\left(z\right)^\ast=\pi_\lambda\left(z\right)^\ast\pi_\lambda\left(i\lambda^{-1}\beta\left(\lambda\right)\right).$$
   Then using the composition rule for $\pi_\lambda$ it follows that
    $e^{-\frac{i}{2}\Re\left(z.\overline{\beta\left(\lambda\right)}\right)}= e^{\frac{i}{2}\Re\left(z.\overline{\beta\left(\lambda\right)}\right)}, ~~~\text{   for all } z\in \C^n.$
    Thus, $\beta(\lambda)=0$ which proves our assertion.
\end{proof}

\noindent{\bf Acknowledgements.} 

\noindent
This work is a part of first author's Ph.D. thesis. He thanks the National Institute of Science Education and Research Bhubaneswar, India, for providing an excellent research facility. RM acknowledges the partial support provided by the research grants  (ANRF/ARGM/2025/000750/MTR).


\smallskip
\noindent{\bf Data Availability Statement.} Not applicable.

\smallskip
\noindent{\bf Code Availability Statement.} Not applicable.


\smallskip

\bibliography{ref12}
\bibliographystyle{amsplain}
\end{document}